\documentclass[12pt]{article}

\author{D.~Zvonkine\thanks{b\^at. 425,
Universit\'e Paris-Sud, 91400 Orsay, FRANCE. \,
{\bf E-mail:} dimitri.zvonkine@ math.u-psud.fr
\hspace{5mm}
The author is partially supported by the RFBR
grant 02-01-22004.}}

\title{Counting ramified coverings and intersection
theory on Hurwitz spaces II \\
{\Large (Local structure of Hurwitz spaces and \\
combinatorial results)}}

\date{\today}

\usepackage{amssymb}
\usepackage{theorem}

\def\A{{\cal A}}
\def\C{{\mathbb C}}
\def\Z{{\mathbb Z}}
\def\Aut{{\rm Aut}}

\def\P{{\mathbb P}}
\def\CP{{\mathbb C} {\rm P}}

\def\dim{{\rm dim}}

\def\qed{{\hfill $\diamond$}}
\def\M{{\overline M}}

\def\A{{\cal A}}
\def\cH{{\cal H}}

\def\cM{{\cal M}}

\def\ocM{{\overline \cM}}

\newtheorem{theorem}{Theorem}
\newtheorem{proposition}{Proposition}[section]
\newtheorem{corollary}[proposition]{Corollary}

\newtheorem{lemma}[proposition]{Lemma}
{\theorembodyfont{\rmfamily}
\newtheorem{conjecture}[proposition]{Conjecture}
\newtheorem{definition}[proposition]{Definition}
\newtheorem{notation}[proposition]{Notation}
\newtheorem{example}[proposition]{Example}
\newtheorem{remark}[proposition]{Remark}
}

\include{epsf}

\begin{document}

\maketitle

\begin{abstract}
The Hurwitz space is a compactification of the space
of rational functions of a given degree. We study
the intersection of various strata of this space with
its boundary. A study of the cohomology ring of the
Hurwitz space then allows us to obtain recurrence
relations for certain numbers of ramified coverings of
a sphere by a sphere with prescribed ramifications.
Generating functions for these numbers
belong to a very particular subalgebra of the algebra of
power series.
\end{abstract}

\section{Introduction}

\subsection{Preliminaries}

This paper is a continuation of the paper~\cite{LZ-I} by
S.~Lando and the author.
Here we use the same general framework and the same notations.

In~\cite{LZ-I} we introduced Hurwitz spaces $\cH_n$
and their projectivizations $\P\cH_n$ (along the lines
of~\cite{ELSV}) and the Lyashko-Looijenga map $LL$. We
also described various strata in Hurwitz spaces 
and proved several relations satisfied by the homology 
classes represented by these strata. 
The degree of the Lyashko-Looijenga
map on each stratum and the corresponding Hurwitz number
were related to the coupling of the stratum
with a power of a particular $2$-cohomology class
$\Psi_n \in H^2(\P\cH_n)$.

A {\em stratum} $\Sigma_{\kappa_1, \dots, \kappa_c} \in \cH_n$ is
described by a set of $c$ partitions of positive
integers $d_1, \dots, d_c$. A generic function lying in this
stratum has $c$ critical values of multiplicities
$d_1, \dots, d_c$. The elements of the partitions correspond
to the multiplicities of the critical points. In the sequel
we often omit in the above 
notation the trivial partitions $\kappa = 1$
corresponding to simple critical points. We denote
by $\mu_{\kappa_1, \dots, \kappa_c}$ the degree of the
$LL$ map on the stratum $\Sigma_{\kappa_1, \dots, \kappa_c}$.

A {\em Hurwitz number} $h_{\kappa_1, \dots, \kappa_c}$ is
the number of ramified coverings with $c$ prescribed
ramification points and ramification types $\kappa_1, 
\dots, \kappa_c$, each covering being counted with weight
$1/|\Aut|$ (where $|\Aut|$ is the number of automorphisms
of the covering).

The precise relation formula obtained in~\cite{LZ-I} 
is the following:
$$
\frac{h_{\kappa_1, \dots, \kappa_c}}
{|\Aut \{\kappa_1, \dots, \kappa_c\}|}
=
\frac{\mu_{\kappa_1, \dots, \kappa_c}}
{n! \, |\Aut \{ d_1, \dots, d_c\}|}
=
\frac{1}{n!} \, 
\left< 
\Psi_n^d, \P \Sigma_{\kappa_1, \dots, \kappa_c}
\right>,
$$
where $d$ is the dimension of the projectivized stratum
$\P\Sigma_{\kappa_1, \dots, \kappa_c}$.

Codimension $1$ strata have special names:
the {\em caustic} $C_n = \Sigma_2$ and the 
{\em Maxwell stratum} $M_n = \Sigma_{1^2}$. 
We also use boundary strata $\Delta_n$ and $\Delta_{p,q}$
(see~\cite{LZ-I}).

\subsection{The results}

In the present paper we use the above results to derive 
recurrence relations (expressed as partial differential
equations on generating functions) on some Hurwitz numbers
corresponding to certain strata in the Hurwitz space.
These relations have a meaning similar to the WDVV equation 
in quantum cohomology: intersecting
a stratum with the boundary of the Hurwitz space, we split
the initial Riemann sphere into several spheres with
smaller numbers of marked points.
These results are formulated in Section~\ref{Sec:recurrence}.

Let $\Sigma$ be a stratum 
(i.e., the closure of the set of rational functions with
prescribed multiplicities of critical points and values)
in the Hurwitz space $\cH_n$.
We first study the intersection of $\Sigma$ with the boundary
of the Hurwitz space (i.e., with the set of stable
rational functions defined on singular nodal curves). 
Our main geometrical result
is that the neighborhood of this intersection locally looks
like a smaller Hurwitz space (Theorem~\ref{Thm:neighborhood}).
In particular, this allows us to find the multiplicity
of the intersection. This investigation of the geometry
of the Hurwitz space is carried out in Section~\ref{Sec:int}.

In Section~\ref{Sec:proofs} we use the cohomological
relations of~\cite{LZ-I} and the results of
Section~\ref{Sec:int} to prove the recurrence relations
of Section~\ref{Sec:recurrence}.

Finally, in Section~\ref{y1y2} we give some explicit
formulas for Hurwitz numbers, resulting from our recurrence relations. 
Here are some examples of these enumerative results.
Let $\kappa_1, \dots, \kappa_c$ be partitions such that
the sum of all their elements equals $2n-2$. 
As above, $h_{\kappa_1, \dots, \kappa_c}$ is the corresponding
Hurwitz number. Let $\alpha$ be the partition $2$
and $\beta$ the partition $3$. They correspond to critical
values whose preimage contains, respectively, a unique
double critical point or a unique triple critical
point. In the notation
below we omit the partition $1$ that corresponds
to a simple ramification point.

\begin{theorem}
\begin{eqnarray}
\nonumber
h_{\alpha, \alpha} & = & 
\frac34  \, (27n^2-137n+180) \,
\frac{n^{n-6} \, (2n-6)!}{(n-3)!}.\\
\nonumber
h_{\alpha, \beta} & = & 
8 \, (6n^2-37n+60) \,
\frac{n^{n-7}\, (2n-7)!}{(n-4)!}.\\
\nonumber
h_{\beta, \beta} & = &
\frac49 \, (256n^3 - 2787n^2 + 10448n - 13440) \,
\frac{n^{n-8}\, (2n-8)!}{(n-4)!}.\\
\nonumber
h_{\alpha, \alpha, \alpha} & = & 
\frac18 \, (729n^3 - 6723n^2 + 21026n - 22680) \,
\frac{n^{n-7}\, (2n-8)!}{(n-4)!}.
\end{eqnarray}
\end{theorem}

We have an algorithm that allows one to find similar formulas
for other cases and a conjecture about the general form
of the answer.

We also prove that some generating series for the numbers
of ramified coverings belong to a very special
subalgebra of the algebra of formal power series and discuss
the properties of this subalgebra. Namely, fix $c$ points
on the sphere and $c$ partitions $\kappa_1, \dots, \kappa_c$
of positive integers $d_1, \dots, d_c$.
Denote by $h_{\kappa_1, \dots, \kappa_c}(n)$ 
the number of the $n$-sheeted coverings
of the sphere by a sphere, ramified over the $c$ points
with ramification types $\kappa_i$ and,
in addition, having $d(n) = 2n-2 - \sum d_i$ fixed  simple ramification
points. Denote by $f_{\kappa_1, \dots, \kappa_c}$ the series
$$
f_{\kappa_1, \dots, \kappa_c}(t) = 
\sum \frac{h_{\kappa_1, \dots, \kappa_c}(n)}{d(n)!} \, t^n.
$$

\begin{theorem} \label{Thm:algebra}
The power series $f_{\kappa_1, \dots, \kappa_c}$ 
belongs to the algebra generated by the power series
$$
\sum_{n \geq 1} \frac{n^{n-1}}{n!} t^n \qquad \mbox{and}
\qquad \sum_{n \geq 1} \frac{n^n}{n!} t^n.
$$
\end{theorem}

We believe that this
subalgebra may play an important role in the combinatorial
theory of ramified coverings of the sphere.

\subsection{Acknowledgments}

The author is grateful to S.~Lando for many useful discussions
and for his remarks on the text. 
J.-M.~Bismut, D.~Panov, and N.~Markarian contributed to 
the work by their questions and helpful remarks.
I also thank V.~Arnold,
A.~Mednykh, A.~Zorich and A.~Zvonkin for their interest.

\section{Recurrence relations for Hurwitz numbers}
\label{Sec:recurrence}

Here we write out recurrence relations
for some Hurwitz numbers. We introduce generating
functions whose coefficients are these Hurwitz numbers.
The recurrence relations
can then be expressed as partial differential equations
with an infinite number of terms for each generating function.
These relations are proved in the two following sections
(Sections~\ref{Sec:int} and~\ref{Sec:proofs}).

Let $\Sigma = \Sigma_{\kappa_1, \dots, \kappa_c}$ be a stratum in
$\cH_n$.
The main idea of the proof is to combine the expression
$$
h_{\kappa_1, \dots, \kappa_c} = \frac1{n!} \, 
|\Aut \{\kappa_1, \dots, \kappa_c\}| \,\,
\left< \P\Sigma \, , \, \Psi_n^{\dim \P\Sigma} \right>
$$
for the Hurwitz numbers (Theorems~3 and~4 in~\cite{LZ-I})
and the expression
$$
\Psi_n = \frac1{2n(n-1)} \,\, \sum_{p+q=n} pq \, \Delta_{p,q}
$$
for the cohomology class $\Psi_n$ (Proposition~4.9 in~\cite{LZ-I}).
From these two identities we deduce:

\begin{corollary}\label{Cor:recurrence}
$$
2n(n-1) \frac{h_{\kappa_1, \dots, \kappa_c}}
{|\Aut \{\kappa_1, \dots, \kappa_c\}|} = 
\frac1{n!} \sum_{p+q=n} pq \; 
\left< (\P\Sigma \cap \Delta_{p,q}) \, , \, 
\Psi_n^{\dim \P\Sigma - 1} \right> \,\, .
$$
\end{corollary}

The right-hand side of the last equality 
can be interpreted as a sum of degrees 
of the $LL$ map on the subvarieties 
of the from $\Sigma \cap \Delta_{p,q}$. A point of this
intersection is a meromorphic function defined on a curve with at
least 2 irreducible components. We will see that it can happen that
the curve actually consists of more than $2$ components,
even for a generic point of the intersection.
Consider a particular irreducible component of the
curve, and suppose that the restrictions of the stable
meromorphic functions to this component are of degree $m$.
Then these restrictions form a stratum in the space
$\cH_m$. Knowing the degrees of the 
$LL$ map on such strata for each component of the curve
(or, equivalently, the corresponding Hurwitz numbers), 
one can find the
degree of $LL$ on the intersection $\Sigma \cap \Delta_{p,q}$.

In order to be able to use this method, we will restrict
our attention to particular classes of strata $\Sigma$, whose
intersection with $\Delta_{p,q}$ is not too difficult
to understand. Moreover, these classes must satisfy the
following property.

\begin{definition}
A class of strata $\Sigma$ is called {\em stable under
intersection with $\Delta$} if, for any
$p,q$, the restrictions of the functions 
$f \in \Sigma \cap \Delta_{p,q}$ to any given
component of the nodal curve, form a stratum $\Sigma'$
that belongs to the same class.
\end{definition}

Here is a list of several classes of strata stable
under intersection with $\Delta$. The classes are
described by properties that should be satisfied
by a generic function of a stratum $\Sigma$ in order
for $\Sigma$ to belong to the class.

\begin{enumerate}

\item All strata.

\item All critical values (and therefore all
critical points) are simple.

\item All critical values are simple except maybe one.

\item There is a critical point with the maximal
possible multiplicity (i.e., multiplicity $n-1$ for
functions of degree $n$).

\item The values at distinct critical points are distinct
(but the critical points can be multiple).

\item The values at distinct critical points are distinct
and the multiplicities of the critical points are equal
to $1$ or $2$.

\end{enumerate}

\subsection{Strata with distinct critical values at
distinct critical points}
\label{simple}

In this subsection we consider ramified coverings
with particular ramification types.
Namely, we consider rational functions $f$ such that
the preimage under $f$ of every critical value of $f$
contains only one critical point.
We introduce a generating function for the corresponding
Hurwitz numbers and write down
a partial differential equation on this generating function.
A closed formula for such Hurwitz numbers is still
unknown, but we have found many new formulas for
particular cases.

\begin{definition}
A stratum $\Sigma$ is called {\em simple} if for
a generic function $f \in \Sigma$, every critical value
of $f$ is attained at a unique (possibly multiple)
critical point.
In other words, the set of partitions $\kappa_1, \dots, \kappa_c$
that describes the stratum $\Sigma$ (see Introduction)
contains only partitions consisting of a single element.
\end{definition}

\begin{example}
The caustic is a simple stratum, while
the Maxwell stratum is not. 
\end{example}

Before proceeding we briefly describe
the intersection of a simple stratum $\Sigma$ with
$\Delta_{p,q}$. Consider a function $f$ in the
open part of $\Sigma$ that tends to a stable meromorphic
function in the intersection $\Sigma \cap \Delta_{p,q}$.
Then, generically, two critical values of $f$ get
glued together. 

The monodromies corresponding to
these critical values are cyclic permutations (i.e.,
have only one cycle of length greater than~$1$), because
each critical value has a unique critical point in
its preimage. When the critical values are glued together
the monodromies are multiplied. The two cycles
in the monodromies have $k \geq 2$ elements in common and
their product is a permutation with $k$ cycles. One can see
that this implies that the rational curve splits into $k$
components on which the meromorphic function is not constant.
If $k \geq 3$, there is one additional component, on
which the function is a constant. The other $k$ components
are glued to this additional one. The cases
$k = 0$ and $k=1$ do not interest us, because one can see
that in these cases the rational curve does not split.

\bigskip

We introduce the following notation for the simple strata.
Let
$$
Y = y_1^{m_1} y_2^{m_2} \dots y_{n-1}^{m_{n-1}}
$$
be a monomial of degree $2n-2$
in variables $y_1, \dots y_{n-1}$, the variable $y_i$
having degree $i$, i.e.,
$$
\sum i \, m_i = 2n-2.
$$
Denote by $\Sigma(Y)$ the simple stratum such that a generic
function $f \in \Sigma(Y)$ has $m_1$ simple critical points,
$m_2$ double critical points, \dots, $m_{n-1}$ critical
points of order $n-1$. Denote by $h(Y)$ the Hurwitz number
corresponding to the stratum $\Sigma$.
Finally, denote by $|\Aut (Y)|$ the number
$$
m_1! \dots m_{n-1}!
$$

Consider the following generating function:
$$
G(t,y_1,y_2, \dots)
=
\sum_{n=1}^{\infty} \sum_{Y}
\frac{h(Y)}{|\Aut (Y)|} \; t^n \, Y.
$$
The second sum is taken over the monomials $Y$ of degree $2n-2$.

The first terms of $G$ are
$$
G = t +
\left(
\frac14 y_1^2
\right) t^2
+
\left(
\frac16 y_1^4 + \frac12 y_1^2y_2 +\frac16 y_2^2
\right) t^3
+
$$
$$
+
\left(
\frac16 y_1^6 + \frac98 y_1^4y_2 + \frac23 y_1^3y_3
+\frac32 y_1^2 y_2^2 + y_1y_2y_3 + \frac16 y_2^3+
\frac18 y_3^2
\right) t^4 + \dots \; .
$$
The first term $t$ is added by convention; it does not
really correspond to a stratum.

\bigskip

We are going to prove that the function $G$ satisfies
a partial differential equation with an infinite number of
terms. In order to describe the terms we introduce the
notion of a decomposition list. To each decomposition list
corresponds a term of the partial differential equation.

\begin{definition}\label{Def:admissible}
A list of non-negative numbers $(a_1, \dots a_k; i,j)$
is called a {\em decomposition list} if
it is possible for two cyclic
permutations with cycles of lengths $i+1$ and $j+1$
to have a product with cycle lengths $a_1+1$, \dots,
$a_k+1$. 
\end{definition}

More precisely, consider two permutations
$\sigma_1$ and $\sigma_2$ that are cycles of lengths
$i+1$ and $j+1$. We suppose that $\sigma_1$ and $\sigma_2$
act on $i+j+2-k$ elements and their cycles have
$k \geq 2$ elements in common. Thus each of the $i+j+2-k$
elements participates in at least one of the cycles.
Suppose that the product $\sigma_1 \sigma_2$ 
has $k$ cycles of lengths $a_1+1, \dots, a_k+1$.
Then $(a_1, \dots a_k; i,j)$ is a decomposition list.
Conversely, for any decomposition list one can find
two permutations satisfying the above description.

This definition is equivalent to the following conditions:
$$
2 \leq k \leq \min (i+1,j+1)
$$
and
$$
a_1 + \dots + a_k = i+j - 2k+2.
$$

Now to each decomposition list we assign three positive
integers (related to each other).

\begin{definition} \label{Def:p(L)q(L)}
Let $L = (a_1, \dots a_k; i,j)$ be a decomposition list. 
Denote by $q(L)$ the number of ways to choose $k$
integers $b_i$, $0 \leq b_i \leq a_i$ whose sum is
equal to $i-k+1$. Denote by $p(L)$ the number
$$
p(L) = (a_1+1) \dots (a_k+1) \, q(L).
$$
Finally, denote by $m(L)$ the number
$m(L) = p(L)/(a_k+1) .$
\end{definition}

The numbers $q(L)$ and $p(L)$ are symmetric with respect
to permutations of the $a_r$ and with respect to the
transposition of $i$ and $j$.
The number $p(L)$ is the number of ways to
decompose a given permutation with cycles of lengths
$a_1+1$, \dots, $a_k+1$ into a product of two cyclic
permutations with cycles of lengths $i+1$ and $j+1$.
The numbers $m(L)$ for various decomposition lists $L$
will play the role of coefficients in the equation on $G$.

Denote by $D_0$ the differential operator 
$D_0 = t \, \partial/\partial t$ and by $D_i$, $i \geq 1$,
the differential operator $D_i = \partial/\partial y_i$.
Instead of writing $D_i f$ we will usually right
$f_i$. Thus $(G_1 G_2)_0$ means 
$t \; \partial/\partial t \;
( \partial G/\partial y_1 \cdot \partial G/ \partial y_2)$.

\begin{theorem}\label{Thmmain2}
The generating function $G$ satisfies the following
partial differential equation.
$$
2t^2 \frac{\partial^2G}{\partial t^2}
=
\sum m(L) 
\; y_i y_j \; (G_{a_1} \dots G_{a_{k-1}})_0 \; 
(G_{a_k})_0,
$$
where the sum is taken over all decomposition
lists $L = (a_1, \dots, a_k;i,j)$.
\end{theorem}

The first terms of the partial differential equation
are:
$$
\begin{array}{l}
2t^2  \frac{\displaystyle \partial^2 G}{\displaystyle
\partial t^2} = \\
\ \\
= \; y_1^2 \; \biggl((G_0)_0 \biggr)^2 
\; + \; 6 \, y_1y_2 \; (G_0)_0 \; (G_1)_0 
\; + \\
\ \\
+ \; 4\, y_2^2 \; (G_0)_0 \; (G_2)_0 \; + \; 4\, y_2^2 
\; \biggl( (G_1)_0 \biggr)^2
\; + \; y_2^2 \; ((G_0)^2)_0 \; (G_0)_0 \; + \\
\ \\
+ \; 8\, y_1y_3 \; 
(G_0)_0 \; (G_2)_0 \; + \; 
4\, y_1y_3 \; \biggl( (G_1)_0 \biggr)^2
\;+ \\
\ \\
+ \;2\, y_2 y_3 \; ((G_0)^2)_0 \; (G_1)_0 \; 
+ \; 8\, y_2 y_3 \; (G_0 G_1)_0 \; (G_0)_0 \; + \\
\ \\
+\; 10\, y_2 y_3 \;(G_0)_0 \; (G_3)_0 \;
+ \;20\, y_2 y_3 \;(G_1)_0 \; (G_2)_0 \;+ \;\dots
\end{array}
$$
Here we have added up the coefficients of similar
terms arising from decomposition lists that differ only in
the order of the numbers $a_1,\dots,a_{k-1}$ or in
the order of $i$ and $j$.

In Section~\ref{y1y2} we give some explicit formulas
for the coefficients of $G$ resulting from 
Theorem~\ref{Thmmain2}. Some of them were known
before, but many are new.

\subsection{Strata with at most one multiple critical value}
\label{primitive}

In this subsection we consider the strata $\Sigma$ such that
a generic function in $\Sigma$ has only one multiple critical
value (and any number of simple critical values). The
degrees of $LL$ on such strata, as well as the
corresponding Hurwitz numbers, were first given
by Hurwitz without a complete proof. Today several
proofs are known: a reconstruction of Hurwitz's proof
by V.~Strehl~\cite{Strehl}, a combinatorial proof by I.~Goulden and
D.~Jackson~\cite{GouJac}, and an algebro-geometric proof by 
T.~Ekedahl, S.~Lando, M.~Shapiro and A.~Vainstein~\cite{ELSV}.
Here we give a new algebro-geometric proof.

More precisely, we will introduce a generating
function $F$ whose coefficients are the above Hurwitz numbers
and prove that it satisfies a partial differential
equation with an infinite number of terms. This
equation allows one to find all the coefficients
of $F$.

\begin{definition}
A stratum $\Sigma$ is called {\em primitive} if a generic
function $f \in \Sigma$ has a unique multiple critical
value (that can be attained at several critical
points of different multiplicities). In other words,
the set of partitions $\kappa_1, \dots, \kappa_c$
that describes the stratum contains at most one
partition different from $1$.
\end{definition}

As above, before proceeding we will briefly describe the intersection
of a primitive stratum $\Sigma$ with $\Delta_{p,q}$. Consider a
function $f$ in the open part of $\Sigma$, that tends to a stable
meromorphic function in $\Sigma \cap \Delta_{p,q}$. Generically,
two critical values of $f$ get glued together. There are
two possibilities. 

First, the two critical values can be
simple. In that case the rational curve splits into two
components. The preimages of the multiple critical
value are distributed (in some way) between these two components.
This situation corresponds to the first term on the
right-hand side of the equality in Theorem~\ref{Thmmain1}.

Second, the multiple critical value can get glued with
a simple critical value. In that case the rational curve
also splits into two components. But now one
of the preimages of the multiple critical value
tends to the intersection point of the two components.
This intersection point is, in general, a critical point for the
functions on both components. Such a situation corresponds
to the second term (the infinite sum) on the right-hand
side of the equality in Theorem~\ref{Thmmain1}.

\bigskip

Primitive strata can be described by monomials in 
an infinite number of variables $x_1, x_2, \dots$.
Let $X = x_1^{m_1} \dots x_n^{m_n}$ be a monomial of
degree $n$, the variable $x_i$ having degree $i$,
i.e, $\sum i \, m_i =n$. 
We denote by $\Sigma(X)$ the primitive stratum
in $\cH_n$ such that the unique multiple critical value
has $m_1$ simple preimages, $m_2$ double preimages,
\dots, $m_n$ preimages of multiplicity $n$. 

To the monomial $X$ we assign the following numbers.
The number $|\Aut (X)|$
is the number $m_1! \dots m_n!\,$. The number
$R(X)$ is the number of simple critical values
of a generic function in the stratum $\Sigma(X)$; it is equal
to $n+p-2$, where $p = \sum m_i$ 
is the number of factors in the monomial $X$.
Finally $h(X)$ is the number of ramified coverings
of the sphere with $R(X)$ simple ramification points
and one ramification point with ramification type
given by the partition $\kappa = 1^{m_2} \dots (n-1)^{m_n}$.
(In other words, $h(X) = h_{\kappa}$.)

Note that, unfortunately, a preimage of multiplicity $k \geq 2$ 
is a critical point of multiplicity $k-1$.

We introduce the following generating function
$$
F(x_1, x_2, \dots) = 
\sum_{n=1}^{\infty} \sum_X 
\frac{n \; h(X)}{R(X)! \; |\Aut (X)|} \; X \, t^n,
$$
where the second sum is taken over the monomials $X$ of
degree $n$.

Hurwitz's formula for the the Hurwitz numbers $h(X)$ 
gives the following expression for $F$:
$$
F = 
\sum_{n=1}^{\infty}
\sum_{p=1}^n \>
\sum_{k_1 + \dots + k_p = n} 
\frac{n^{p-2}}{p!}\;
\prod_{i=1}^p \frac{k_i^{k_i}}{k_i!} \,x_{k_i} \; t^n . 
$$
The first terms of $F$ are
$$
F = tx_1 + 
t^2\left(
\frac12 x_1^2+x_2
\right)
+
t^3\left(
\frac12 x_1^3 + 2 x_1x_2 + \frac32 x_3
\right)
+
$$
$$
+t^4\left(
\frac23 x_1^4 + 4 x_1^2x_2 + \frac92 x_1x_3 + 2 x_2^2 + \frac83 x_4
\right) + \dots 
$$
The very first term $t x_1$ is added by convention.
It does not really correspond to a stratum.

\begin{theorem}\label{Thmmain1}
The function $F$ satisfies the following partial differential
equation
$$
2 t^2 \frac{\partial}{\partial t} 
\left( \frac{F}{t} \right)
=
\left( t \frac{\partial F}{\partial t} \right)^2
+ \sum_{i,j \geq 1}
(i+j) x_{i+j} 
\frac{\partial F}{\partial x_i}
\frac{\partial F}{\partial x_j}
$$
\end{theorem}

This equation and the initial condition $F = tx_1 + \dots$
determine all the coefficients of $F$. The proof of the
theorem is completely independent of the Hurwitz formula.

\begin{remark}
Goulden and Jackson \cite{GouJac} obtained by purely
combinatorial methods a similar partial
differential equation for a generating function that
enumerates the same numbers. The idea of their method is
as follows. The Hurwitz number $h(X)$ 
with $X = x_{k_1} \dots x_{k_p}$
is equal (up to a simple combinatorial
factor) to the number of shortest possible transitive
factorizations 
into transpositions of a given permutation $\sigma \in S_n$
whose cycles have lengths $k_1, \dots, k_p$. Multiplying
the last transposition of the factorization by $\sigma$,
we obtain a new permutation $\sigma'$, which is again
factorized into transpositions, but not necessarily 
transitively (i.e., the group generated by the transpositions
may not be transitive). Considering various possible
cases one can express the number of shortest transitive
factorizations for the initial stratum in terms of
similar numbers for several other strata, which is
then encoded in a partial differential equation on
the generating function. From the algebro-geometric
point of view, the case when the factorization of 
$\sigma'$ is not transitive corresponds to the
situation when the curve of definition of the 
meromorphic function $f$ splits into 2 components.
The case when the factorization of $\sigma'$ remains
transitive corresponds to other possible degenerations.
Thus our equation is different from that of
Goulden and Jackson in that we have managed to get
rid of all the degenerations except the splitting of
the rational curve. The price to pay is that our
equation contains nontrivial coefficients coming
from the algebraic geometry (they are multiplicities
of intersections of some strata). Therefore we think
that our equation cannot be easily obtained by purely
combinatorial methods.
\end{remark}

\section{Studying the intersection of a stratum $\Sigma$
with $\Delta_{p,q}$}
\label{Sec:int}

The proofs of Theorems~\ref{Thmmain2} and~\ref{Thmmain1}
are obtained by studying the intersection of primitive
and simple strata with $\Delta_{p,q}$. (Recall that
$\Delta_{p,q}$ is the closure of the set of stable
functions defined on curves with two irreducible
components, the functions being of degree $p$ on
one component and $q$ on the other one.)
This intersection usually consists of a very large number
of irreducible components. The multiplicities of the
intersection can be different for different irreducible
components. Below the word {\em component} means a union
of irreducible components.

For a precise definition of the 
strata $\Delta_{p,q}$ and of the cohomology class $\Psi_n$
see in~\cite{LZ-I}, Definitions~4.8 and 2.11.

\subsection{The $k$-eared strata}
\label{Sub:k-eared}

Here we introduce the notion of $k$-eared stratum and
prove that every component of the intersection of a
simple or primitive stratum $\Sigma$ with $\Delta_{p,q}$ is
contained in a $k$-eared stratum.
Let $\Sigma$ be any stratum in $\cH_n$. Denote by $d$ the 
dimension of $\P\Sigma$; thus a generic function $f$ in $\Sigma$
has $d+1$ critical values.

\begin{lemma} \label{Lem:numcritval}
Let $\Sigma'$
be a component of the intersection $\Sigma \cap \Delta_{p,q}$. Suppose
that a generic function in $\Sigma'$ has less than $d$ critical
values. Then $\left<\P\Sigma', \Psi_n^{d-1} \right> = 0$.
\end{lemma}

\begin{remark}
Since our goal is to find the sum of 
$\left< \P\Sigma', \Psi_n^{d-1} \right>$
over all the components of the intersection $\Sigma \cap \Delta_{p,q}$,
what the lemma actually says is that we can 
discard the components where 
more than two critical values get glued together.

Actually, there are no components satisfying
the conditions of the lemma, but this fact is
much more difficult to prove than the lemma itself
and the lemma is sufficient for our purposes.
\end{remark}

\paragraph{Proof of the lemma.} 
The image of $\Sigma'$ under $LL$ has dimension
at most $d-1$ (i.e., less than the dimension of $\Sigma'$
itself). On the other hand, the 2-cohomology class $\Psi_n$ 
is the pull-back of a 2-cohomology class on the image of $LL$
(see the proof of Theorem~4 in~\cite{LZ-I}).
Therefore the number $\left< \P\Sigma', \Psi_n^{d-1} \right>$
can be computed in the image of $LL$, where it is equal to 0
for dimension reasons.
\qed

\bigskip

Now we will consider a particular kind of subvarieties of
$\cH_n$, called the $k$-eared strata, 
and study their neighborhoods in $\cH_n$. As we will soon
see, all the intersections of $\Delta_{p,q}$ with simple and
primitive strata lie in such subvarieties. Therefore, studying
their neighborhoods will help us to find the multiplicities
of the intersections of $\Delta_{p,q}$ with simple and
primitive strata.

A $2$-eared stratum is consists of stable meromorphic functions
defined on a $2$-component curve, such that the restrictions
of the functions to each component belong to some given strata,
while the common point of the two components is 
a critical point of given multiplicity
on each component of the curve.

A $k$-eared stratum for $k \geq 3$ 
is consists of stable meromorphic functions
defined on a $(k+1)$-component curve. One of the components is
{\em central}. The meromorphic function is constant on this
component. All the other $k$ components are called 
{\em peripheral}. They do not intersect each other, but
each of them intersects the central component
at one point. The restriction of the meromorphic function to
a peripheral component belongs to a given stratum. The point
of intersection of a peripheral component with the central
component is a critical point of given multiplicity.

In order to give a precise definition of a $k$-eared
stratum, let us first define a stratum with a distinguished
point. Recall that a stratum is described by an unordered
set of partitions (see Definition~3.3 in~\cite{LZ-I}),
each partition corresponding to a critical value
of the stable function, and an element $i$ in a partition 
corresponding
to a critical point of multiplicity $i$. We can similarly
describe a stratum of stable functions with a distinguished
point on the curve $S$. In order to encode the necessary information
about the distinguished point we must make the following choices. 
If we want the value of the function at the distinguished point
to be a critical value, we chose one of the partitions $\kappa$
in the set that describes the stratum. Otherwise, if we
want the value at the distinguished point to be an
ordinary value, we do not choose anything. Suppose we have
chosen a partition $\kappa$. Now we must decide whether the distinguished
point will be a critical point or not. In the first case
we choose one of the elements $i$ of the partition; in the
second case we choose nothing. Thus our choice encodes one
of the following possibilities: either the distinguished
point is a critical point of multiplicity $i$ with critical
value of type $\kappa$; or the distinguished
point is a noncritical point, but the value of the function
at the distinguished point is a critical
value of type $\kappa$; or the value of the function at the
distinguished point is noncritical.

Now a stratum of functions with a distinguished point is (the
closure of) the set of all stable functions defined on nodal curves $S$
with one distinguished point, the distinguished point 
satisfying the above conditions.

Let $k \geq 2$ be a natural number. 
Divide $n$ into a sum of $k$
natural numbers $n = n_1 + \dots + n_k$, $n_i \geq 1$.
Divide the set of poles $\{1, \dots, n\}$ into $k$ parts
with $n_1, \dots, n_k$ elements. Choose $k$ strata
$\Sigma_1, \dots, \Sigma_k$ of stable functions with a distinguished
point in the spaces $\cH_{n_1}, \dots, \cH_{n_k}$.

\begin{definition}
Let $k=2$. Consider the set of all stable meromorphic
functions defined on 2-component curves, such that the
restrictions of the function to the components belong to
open parts of the strata $\Sigma_1$, $\Sigma_2$ and the intersection
point of the components is the distinguished point for both
components. This set is called an {\em open 2-eared stratum}.
Its closure in $\cH_n$ is called a {\em 2-eared stratum}.

Let $k \geq 3$. Consider the set of all stable meromorphic
functions defined on $(k+1)$-component curves with one
central and $k$ peripheral components, such that the
restriction of the function to the central component is
a constant, while its restrictions to the peripheral
components belong to
open parts of the strata $\Sigma_1$, \dots, $\Sigma_k$, 
while the intersection 
point of each peripheral component with the central one
is the distinguished point on the peripheral component. 
This set is called an {\em open k-eared stratum}.
Its closure in $\cH_n$ is called a {\em k-eared stratum}.
\end{definition}

\begin{proposition}\label{Prop:inter}
Let $\Sigma$ be a simple or a primitive stratum. Consider
an irreducible component $\Sigma'$ of the intersection 
$\Sigma \cap \Delta_{p,q}$ such that
$$
\left< \P\Sigma', \Psi_n^{\dim \, \P\Sigma'} \right> \not= 0.
$$
Then $\Sigma'$
is contained in a $k$-eared stratum for some $k$. Moreover,
we can, in a natural way, assign a decomposition list $L$ 
{\rm(}Definition~\ref{Def:admissible}{\rm)} to the component 
$\Sigma'$.
\end{proposition}

\paragraph{Proof.} According to Lemma~\ref{Lem:numcritval}
a generic stable function $(S', f') \in \Sigma'$ has one critical
value less than a generic function $(S,f) \in \Sigma$. 
(The function $f'$ is not the derivative of $f$, just
another stable function.)
Therefore as $(S,f)$ approaches  $(S', f')$, two of its critical values
are glued together, and their monodromies are multiplied.

First suppose $\Sigma$ is a primitive stratum. Then $f$ has only
one critical value that is not simple. Therefore, among the
two critical values that get glued together, at least
one is simple, i.e., the corresponding monodromy is a transposition.
When we multiply a transposition by another permutation $\sigma$,
two cases are possible. Either the two permuted elements
of the transposition belong to two different cycles of
$\sigma$, in which case these cycles are merged into one
cycle in the product. Or the two permuted elements
of the transposition belong to the same cycle of
$\sigma$, in which case this cycle splits into two
cycles in the product.

Since $(S', f')$ lies in $\Delta_{p,q}$, we know that the
monodromies of $f'$ do not act transitively on the
set of poles, but have at least two transitivity components
corresponding to the irreducible components of $S'$. It follows
that only the second of the above two cases is possible. It
is easy to see that the set of monodromies of $f'$ has
exactly two transitivity components. Therefore the stable
function $(S', f')$ belongs to a $2$-eared
stratum.

We assign to $\Sigma'$ the decomposition list 
$L = (a_1, a_2; i,1)$,
where $i+1$ is the length of the cycle of $\sigma$ that splits
into two, and $a_1+1$, $a_2+1$ are the lengths of the two
cycles obtained by the splitting. (The last number $1$ means
that the transposition from the above discussion is
a cycle of length $1+1 =2$.)

Now suppose that $\Sigma$ is a simple stratum. Then the monodromies
of the critical points of $f$ that get glued together are
cycles. When two cyclic permutations are multiplied, several
cases are possible. For simplicity,
let us momentarily forget about
the elements of the permutation that do not
belong to either of the two cycles. 
(i) If the cycles do not have common
elements, their product is just a permutation with two cycles
of the same lengths as the cycles that we multiply. (ii)
If the cycles have one common element, then their product
is just one cycle. (iii) If the cycles have $k \geq 2$
elements in common, then their product is a permutation
with $k$ cycles.

Again, we know that the set of monodromies of $f'$
has at least two transitivity components. It follows that
only case (iii) is possible, and that the actual number
of transitivity components is $k$. Thus the curve $S'$
has $k$ components on which $f'$ is not constant. It
remains to check that for $k \geq 3$ it has only one component
on which $f'$ is constant. This follows from the fact
that the elements of two cycles form a unique transitivity
component under the action of these two cycles, provided
the cycles have at least one common element. Thus we have
proved that $(S',f')$ belongs to a $k$-eared stratum. 

We assign to the component $\Sigma'$ the decomposition
list $L = (a_1, \dots, a_k; i,j)$, where $i, j$ are
the multiplicities of the two critical values that are
glued together, and $a_1+1, \dots a_k+1$ are the lengths
of cycles in the product of their monodromies.
\qed

\subsection{Computing the multiplicity of intersection}
\label{Sub:intersection}

Here we study the local geometry of $k$-eared strata, which
allows us to compute the multiplicity of the intersection of
a simple or a primitive stratum $\Sigma$ with $\Delta_{p,q}$.

\begin{proposition}\label{M_k}
Let $k \geq 3$. Consider a $k$-eared stratum $\Sigma$ and a generic
point in its image under the $LL$ map. The preimage of this
point in $\Sigma$ is isomorphic to a finite number of copies of the 
compactified moduli space $\ocM_k$.
\end{proposition}

\paragraph{Proof.} The restriction of the stable meromorphic
function to any peripheral component is determined, up to a
finite number of possibilities, by its image under $LL$.
Its restriction to the central component is a constant,
and the constant is also determined by the image under $LL$.
The only thing that can vary is the disposition on the
central component of the $k$ points at which it intersects 
the peripheral components. These dispositions form the
space $\ocM_k$. Thus each point of a $k$-eared stratum
is naturally contained in a subvariety isomorphic to
$\ocM_k$. Note that the same is true for a projectivized
$k$-eared stratum $\P\Sigma$. \qed

\bigskip

\begin{notation}\label{notation}
We introduce the following notation that the reader must
bear in mind in the sequel:
$$
\ocM_k \equiv A_k \subset B_k \supset \Sigma' 
\subset \Sigma \cap \Delta_{p,q} \subset \cH_n.
$$
Here $\Sigma$ is a primitive or a simple stratum, $\Sigma'$
a component of its intersection with $\Delta_{p,q}$. Further,
$B_k$ is the $k$-eared stratum that contains 
$\Sigma'$ (according
to Proposition~\ref{Prop:inter}), and $A_k$ is a subvariety of
$B_k$ (as in Proposition~\ref{M_k}) containing a point of $\Sigma'$
and isomorphic to the moduli space $\ocM_k$. The image of $A_k$
under the $LL$ map is a point.
In general we will assume that the subvariety $A_k$ is chosen
generically inside $B_k$.
\end{notation}

Now we study the neighborhood of $A_k$ in $\cH_n$.
To do that, we must introduce, following~\cite{ELSV} 
a new kind of Hurwitz spaces,
more general than the one we have used up to now. The
neighborhood of a $k$-eared stratum will be very similar
to such Hurwitz spaces.

Consider a list of nonnegative integers $a_1, \dots, a_k$.
Let $\sigma$ be a permutation of $k + \sum a_r$ elements
with cycle lengths $a_1+1$, \dots, $a_k+1$.

\begin{definition}
We call an {\em indexed Hurwitz space} $\cH(a_1, \dots, a_k)$
the space of all stable meromorphic functions $(S,f)$, 
considered up to an {\em additive} constant, with $k$ poles
of orders $a_1+1, \dots a_k+1$.

We call a {\em decorated Hurwitz space} $\cH_{\sigma}$
the space of the same stable functions as above
endowed with the following additional information. The
preimages under $f$ of the half-line $[A, +\infty]$
are numbered. (Here $A$ is a sufficiently large real number such 
that there are no critical values of $f$ on the
semi-closed interval $[A, +\infty)$.) The monodromy around
$\infty$ (i.e., the permutation of the numbered preimages
obtained by going around $\infty$ in the counterclockwise
direction) is equal to $\sigma^{-1}$.
\end{definition}

It is easy to see that an indexed Hurwitz space
$\cH(a_1, \dots, a_k)$ is the quotient of a decorated
Hurwitz space $\cH_{\sigma}$ by the group
$$
\Z/(a_1+1)\Z \oplus \dots \oplus \Z/(a_k+1) \Z
$$
acting by cyclic renumberings of the preimages
of $[A, +\infty]$ around each pole.

Indexed Hurwitz spaces were studied in detail
in~\cite{ELSV} (for curves of all genera). 
It was shown there that for $k \geq 3$ 
such Hurwitz spaces are cones over $\ocM_k$; in
particular, they are not smooth. Decorated Hurwitz
spaces are smooth and can be seen to be fiber bundles
over $\ocM_k$, as we now briefly describe. 

Indeed, suppose we have fixed the positions
of $k$ points on $\CP^1$. Then the principle part of the stable
meromorphic function at the neighborhood of a pole of
order $a+1$ can be written as
$$
\left( \frac{u}{z} \right)^{a+1} + 
v_1 \left( \frac{u}{z} \right)^a + \dots + v_a \frac{u}{z},
$$
where $z$ is the variable and $(u, v_1, \dots, v_a)$
is a set of parameters. Multiplying $u$ by an $(a+1)$-th
root of unity $\varepsilon$ and each $v_i$ by 
$\varepsilon^i$ amounts to a renumbering of the
preimages of the interval $[A, +\infty]$.
Putting together the above sets of parameters
for all the poles, we obtain a description of a fiber
of the projection from $\cH_{\sigma}$ to $\ocM_k$.
When $u=0$, the point $z=0$ is no longer a pole
of the function. In this case we must glue a new
component to the curve $S$ at the point $z=0$.
The coordinate on this new component will be $w= u/z$.
Thus on this new component a pole of order $a+1$ will
be preserved.

For indexed Hurwitz spaces, the parameters $u$ and $v_i$
are defined only up to multiplication by some powers
of an $(a+1)$-th root of unity. More precisely, 
assign to $u$ the weight $1$
and to $v_i$ the weight $i$. Consider
the algebra of polynomials in the variables $u$ and $v_i$,
such that the weight of each of their monomials is
a multiple of $a+1$. Then the algebra of
functions on a fiber of $\cH(a_1, \dots, a_k)$ over
$\ocM_k$ is the tensor product of such algebras for
$a_1, \dots, a_k$.

\bigskip

Consider a $k$-eared stratum $B_k$ in $\cH_n$. Let $(S,f)$ 
be a generic point of this stratum. Let $w$ be the value
of $f$ on the central component of $S$ and $w'$
a point close to $w$. Let us number those
preimages of $w'$ that lie close to the central component
of $S$. Denote by $\sigma$ the permutation that these
preimages undergo as $w'$ goes around $w$.
Finally, consider the subvariety $A_k$
obtained from the initial generic point of the
$k$-eared stratum by changing the positions of the intersection
points on the central component without touching the
peripheral components.

\begin{theorem}\label{Thm:neighborhood}
For $k \geq 3$, 
the neighborhood of $A_k$ in $\cH_n$ is isomorphic to
the neighborhood of $\ocM_k \times \{0\}$ in
$\cH_{\sigma} \times \C^d$ for some positive integer $d$.
\end{theorem}

\paragraph{Proof.} As in the proof of Proposition~\ref{Prop:inter}
denote by $(S',f')$ a stable function in $A_k$
and by $(S,f)$ a stable function in $\cH_n$ close to it.

The function $(S',f') \in A_k$ is almost a ramified
covering of $\CP^1$ by the curve $S'$. The only difference
with a ramified covering is that the central component of $S'$
is sent to a point by $f'$.
The idea of the proof is that any deformation
of $(S',f')$ can be regarded
as a set of independent deformations of the
ramification structure at the neighborhoods
of the critical values of $f'$. Deforming
the special critical value (the value of $f'$ on the
central component) looks like deforming a point of $\ocM_k$
inside $\cH_{\sigma}$. Deforming any other critical
value looks simply like moving off the origin in a
complex vector space of some dimension.

Consider a point $(S, f)$ in $\cH_n$
close to $(S',f')$. Then, if $W$ is a critical value of 
the function $f'$ whose monodromy is a permutation
$\rho$, the function $f$ will have several critical values
close to $W$, whose product of monodromies is equal
to $\rho$. 

First suppose $W$ is an ordinary critical
value (not the value of $f'$ on the central component
of $S'$). 

The simplest case is when there is only one
critical point in the preimage of $W$, in other words,
the monodromy $\rho$ has only one cycle of length 
$p > 1$. Then at the neighborhood of the critical point
in question one can choose a local coordinate $z$
such that the function $f'$ locally equals $z^p+W$. Now 
a deformation of $f'$ ``at the neighborhood of the
critical value $W$'' can be obtained in the following 
way. Consider a small disc $D_W$ surrounding $W$. Cut out
from the curve $S'$ the connected component of $(f')^{-1}(D_W)$
containing the critical point $z=0$. 
Instead, glue into the hole thus obtained
the preimage of $D_W$ under the function
$z^p + W + (\lambda_1 z^{p-1} + \dots + \lambda_p)$
with sufficiently small complex numbers $\lambda_i$.
Thus we obtain a new curve $S$ with a new stable
meromorphic function $f$ on it. The set of such deformations
is parameterized by the $\lambda_i$s and the initial
function $(S',f')$ corresponds to 
$\lambda_1 = \dots = \lambda_p = 0$.

The complex structure along
the gluing contour on the new surface $S$
is uniquely reconstructed in such a way
that the function $f$ becomes meromorphic.
Alternatively, we can cut out from $S'$
the preimage of a disc slightly smaller than $D_W$
and glue back the preimage of a disc slightly
bigger than $D_W$. The gluing then identifies
two preimages of an open annulus and the complex
structure extends automatically.

Now suppose that there are several critical points
whose image is equal to $W$. This case is actually
just as simple, because we can perform the above
operation of cutting off a small disc and replacing
it by another one independently for all the critical points.

It remains to consider the case when $W=w$ is the value
of $f'$ on the central component of $S'$. First of all,
the preimage of $w$ can contain critical points that do not
lie on the central component of $S'$. These can be deformed
by the same procedure as above. 

Now let us see how the
neighborhood of the central component of $S'$ can be deformed.
First, as above, cut off from the curve $S'$ the connected
component of $(f')^{-1}(D_w)$ containing the central component.
Now consider a stable
meromorphic function $g$ in $\cH_{\sigma}$ close to the base $\ocM_k$; 
this means that the critical values of $g$ are close to each 
other. More precisely, we suppose that by adding an appropriate
constant to $g$ we can make all its critical values
fit inside the disc $D_w$. Then $g$ determines a ramified
covering of the disc $D_w$ with monodromy $\sigma$. Therefore
we can glue the preimage $g^{-1}(D_w)$ instead of $(f')^{-1}(D_w)$
into the curve $S'$. We will obtain a new curve $S$ with a
stable meromorphic function $f$ on it, close to $(S',f')$.

This finishes the proof. \qed

\bigskip

Having described the neighborhood of $A_k$, we describe
the intersections of this neighborhood with the stratum
$\Sigma$ and with $\Delta_{p,q}$. These intersections look
like the varieties $s_L$ and $\delta_r$ below.

Consider an indexed Hurwitz space $\cH(a_1,\dots,a_k)$
and a decorated Hurwitz space $\cH_{\sigma}$, where $\sigma$
is a permutation with cycle lengths $a_1+1, \dots, a_k+1$.
Let $L = (a_1, \dots, a_k; i,j)$
be a decomposition list. 

\begin{definition}
Denote by $s_L \subset \cH_{\sigma}$
the stratum consisting of stable meromorphic functions 
with exactly two critical points of multiplicities
$i$ and $j$. (The corresponding critical values are
automatically distinct.) Denote by $\overline{s}_L \subset
\cH(a_1, \dots, a_k)$ the similarly defined stratum in
the indexed Hurwitz space. (In other words, $\overline{s}_L$
is the quotient of $s_L$ by the action of the 
group $\bigoplus \Z/(a_r+1) \Z$.)

Denote by $\delta_r \subset \cH_{\sigma}$ the stratum
consisting of the stable meromorphic functions defined
on 2-component curves, the $r$th pole being on the
first component,
and all the other poles on the second component.
Denote by $\overline \delta_r \in \cH(a_1, \dots, a_k)$ 
the similarly defined stratum in the indexed Hurwitz space.
(Again $\overline \delta_r$
is the quotient of $\delta_r$ by the action of the
group $\bigoplus \Z/(a_r+1) \Z$.
\end{definition}

Let $(S,f)$ be a generic point of the intersection $\Sigma \cap
\Delta_{p,q}$, where $\Sigma$ is a simple or a primitive
stratum. The point $(S,f)$ is contained in the variety $A_k$ (see
Notation~\ref{notation}). Denote by $L$ the decomposition
list corresponding to the point $(S,f)$ 
(see Proposition~\ref{Prop:inter}).

\begin{proposition} Let $k \geq 3$.

(i) The intersection of the neighborhood of $A_k$ with the stratum
$\Sigma$ is isomorphic to $s(L) \times \C^{d_s}$ for some integer
$d_s$. 

(ii) The intersection of the neighborhood of $A_k$ with $\Delta_{p,q}$ 
is isomorphic to $\bigcup \delta_r \times \C^{d_{\delta}}$ 
for some integers $d_r$, where the union is taken over the numbers
$r \in \{1, \dots, k \}$ such that the restriction of $f$ on the
$r$th peripheral component of $S$ is of degree~$p$ or of degree~$q$.
\end{proposition}

\paragraph{Proof.}

(i) This is a consequence of the proof of 
Theorem~\ref{Thm:neighborhood}. Indeed, to move
off $A_k$ into $\Sigma$ we must deform the critical values
of $f$ in the following way. For an ordinary critical value
(i.e., not the value of $f$ on the central component of $S$):
just move the critical value itself, preserving the
multiplicity of the critical point in its preimage. This
corresponds to moving off the origin in some vector space
$\C^{d_s}$. For the value of $f$ on the central component
of $S$: this value should be decomposed into two critical
values, $k$ cycles of lengths $a_1+1$, \dots, $a_k+1$
being represented as a product of two cycles of lengths
$i+1$ and $j+1$. This is precisely the description of
the stratum $s_L$.

(ii) This, again, follows from the proof of 
Theorem~\ref{Thm:neighborhood}. To move off 
$A_k$ into $\Sigma$ we must deform the critical values
of $f$ in the following way. For an ordinary critical value
(i.e., not the value of $f$ on the central component of $S$):
deform it in any possible way. 
For the value of $f$ on the central component
of $S$: it should be resolved in such a way that only
the $r$th peripheral component (for one of the
numbers $r$ described in the Proposition) of the curve $S$ remains
a separate component, while all the other components
merge together. This is precisely the description
of $\delta_r$. \qed

\bigskip

Now we should compute the multiplicity of the
intersection of $s(L)$ and $\delta_r$. Actually, this
is easier to do for $\overline s(L)$ and $\overline \delta_r$.
We start with a simple geometrical description of
the stratum $\overline s(L)$.

Recall that $p(L)$ is the number
of ways to decompose $\sigma$ into a product of two cycles
with lengths $i+1$ and $j+1$, while
$$
q(L) = p(L)/(a_1+1) \dots (a_k+1)
$$
(see Definition~\ref{Def:p(L)q(L)}).

\begin{proposition}\label{Prop:s_L}
For $k \geq 3$ 
the stratum $\overline s_L \subset \cH(a_1, \dots, a_k)$ 
is a union of $q(L)$
straight lines, each line being contained in a fiber of the
bundle $\cH(a_1, \dots, a_k) \rightarrow \ocM_k$.
\end{proposition}

\paragraph{Proof.} This is almost obvious.
Up to transformations $f \mapsto af+b$, 
there are exactly $p(L)$ stable meromorphic functions
with poles of multiplicity $a_1+1$, \dots $a_k+1$,
with two critical points of multiplicities $i$ and $j$,
with numbered preimages of $[A, +\infty]$ and with
a given monodromy $\sigma$. Forgetting the
numbering, gives us $q(L)$ functions, still up
to transformations $f \mapsto af+b$.
Since a point of $\cH_{\sigma}$ is a meromorphic function
up to an additive constant, the assertion follows. \qed

\begin{proposition} \label{Prop:deltacaps(L)k}
For $k \geq 3$, the multiplicity of intersection of the stratum
$\overline s(L)$ and $\overline \delta_r$ in
$\cH(a_1, \dots, a_k)$ equals
$$
\frac1{a_r+1}\, q(L).
$$

The multiplicity of intersection of the stratum
$s(L)$ and $\delta_k$ (here we have chosen $r=k$) in
$\cH_{\sigma}$ equals $m(L)$ {\rm(}see 
Definition~\ref{Def:p(L)q(L)}{\rm)}.
\end{proposition}

\begin{remark}
The first of the two multiplicities is not necessarily
an integer. This is due to the fact that $\cH(a_1, \dots, a_k)$
is not a smooth variety and $\overline s(L)$ and $\overline \delta_r$
are {\em Cartier divisors},
that is, algebraic subvarieties that are not necessarily
locally equal
to zero sets of functions. Here is the simplest example
of such situation.

Consider the cone $z^2 = x^2+y^2$ in $\C^3$. Let $l_1$ and
$l_2$ be two straight lines belonging to this cone. We claim
that the multiplicity of their intersection is equal to $1/2$.
Indeed, consider the plane tangent to the cone along $l_1$.
Then its intersection with the cone is $2 l_1$. Similarly,
the intersection of the cone with the plane tangent to it
along $l_2$ equals $2l_2$. Now we can move slightly both
planes, and the intersection of the three subvarieties:
the two shifted planes and the
cone will consist of two points. Thus $2 l_1 \cap 2 l_2 = 2$,
so $l_1 \cap l_2 = 1/2$.
\end{remark}

\paragraph{Proof of Proposition~\ref{Prop:deltacaps(L)k}.}
First let $k \geq 3$.

The stratum $\overline \delta_r$ intersects each
fiber of the bundle $\cH(a_1, \dots, a_k)$ transversally,
while the stratum $\overline s_L$ is a union of straight
lines, each of which is contained in one of these fibers.
Therefore we can restrict our attention to one fiber and
find the multiplicity of intersection inside it.

The principle part of the $r$th pole of $f$ can be
written as
$$
\left( \frac{u}{z} \right)^{a_r+1} + 
v_1 \left( \frac{u}{z} \right)^{a_r} + \dots + v_{a_r} \frac{u}{z}.
$$

First consider the space $\C^{a_r+1}$ with coordinates
$u,v_1, \dots, v_{a_r}$ endowed with the action of $\C^*$
$\lambda: u \mapsto \lambda u, \, v_i \mapsto \lambda^i v_i$.
In this space we consider the hyperplane $u=0$ and an
orbit of the $\C^*$ action that is not contained in
this hyperplane. The multiplicity of the intersection 
between the two (at the origin) is equal to $1$.

Now we factor this space by the action of $\Z/(a_r+1) \Z$
(seen as the subgroup of roots of unity of order $a_r+1$
in $\C^*$). The image of the hyperplane is the
subvariety $u^{a_r+1} =0$. The image of the $\C^*$
orbit is a straight line. The multiplicity of the intersection 
between the two equals $1/(a_r+1)$.

It follows that the intersection of a straight line in a fiber
of the bundle $\cH(a_1, \dots, a_k) \rightarrow \ocM_k$
with the stratum $\overline \delta_r$ has multiplicity $1/(a_r+1)$.
Since the stratum $\overline s(L)$ consists of $q(L)$ such lines, 
the first assertion of the proposition follows. The multiplicity
of intersection of $s(L)$ and $\delta$ is equal to
$(a_1+1) \cdot \dots \cdot (a_k+1)$ times the multiplicity
of intersection of $\overline s(L)$ and $\overline \delta_k$, 
whence the second assertion.
\qed

\bigskip

\begin{figure}[h]
\begin{center}
\
\epsfbox{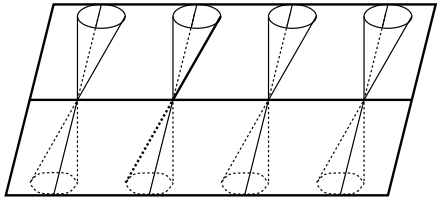}

\caption{\label{Fig:sdelta} The horizontal line is
$\ocM_k$. The space $\cH(a_1, \dots, a_k)$ is a cone
over $\ocM_k$; we have shown several fibers of this cone.
The plane is the stratum $\overline \delta_r$. The bold
oblique line is one of the lines composing $\overline s(L)$.}
\end{center}
\end{figure}

In Figure~\ref{Fig:sdelta} we have symbolically shown the space 
$\cH(a_1, \dots, a_k)$ as a cone over $\M_k$ and the strata
$s(L)$ and $\delta_r$ inside it.

Now we consider the case $k=2$. For $k=2$ the stratum $s(L)$
in $\cH_{\sigma}$ is defined exactly as before, while $\delta$
is now just the set of stable functions defined on $2$-component
curves. One can see that if 
$L = (a_1, a_2; i,j)$ is a decomposition list,
then $q(L) = \min(i,j, a_1+1, a_2+1)$.

\begin{proposition} \label{Prop:deltacaps(L)2}
For $k = 2$, the multiplicity of intersection of the stratum
$\overline s(L)$ and $\overline \delta$ in
$\cH(a_1, \dots, a_k)$ equals
$$
\frac{a_1+a_2+2}{(a_1+1)(a_2+2)} \, q(L).
$$

The multiplicity of intersection of the stratum
$s(L)$ and $\delta$ in $\cH_{\sigma}$ equals $(a_1+a_2+2) q(L)$.
\end{proposition}

\begin{remark}\label{Rem:2andk}
The reader may have noticed that the case $k=2$ has been
excluded from all the previous propositions. Indeed,
since the moduli space $\ocM_2$ does not exist, the
geometry of the intersection of $\Sigma$ with $\Delta_{p,q}$
in the neighborhood of a $k$-eared stratum is different
for $k=2$ and $k \geq3$. However, not surprisingly,
the formula for the multiplicity of the intersection
happens to be the same. Indeed, in the case $k=2$,
the same stratum $\delta$ corresponds to two different
decomposition lists: $(a_1, a_2; i,j)$ and $(a_2, a_1; i,j)$.
(Recall that for $k \geq 3$, the points $(S,f)$ of $\delta_r$ are
stable functions, such that the $r$th pole of $f$ lies on a
separate component of $S$. But for $k=2$, saying that the
second pole lies on a separate component is the same as
saying that the first pole lies on a separate component.)
The corresponding numbers $m(L)$ (Definition~\ref{Def:p(L)q(L)})
are equal to $(a_1+1) q(L)$ and $(a_2+1) q(L)$. Adding them,
we obtain $(a_1+a_2+2) q(L) = (i+j) q(L)$, as in the
assertion of the proposition.
\end{remark}

\paragraph{Proof of Proposition~\ref{Prop:deltacaps(L)2}.}
Let us consider the space $\widehat \cH_{\sigma}$
of stable meromorphic functions $(S,f)$
with two poles of multiplicities $a_1+1$ and $a_2+1$
and with {\em one distinguished point on the curve S}.
This space can be treated as in the previous
proposition (for $k \geq 3$), because it is a fiber
bundle over the moduli space $\ocM_3$, i.e., over a point.
There is a forgetful map from $\widehat \cH_{\sigma}$ to
$\cH_{\sigma}$,
that forgets the distinguished point and, if there appeared 
a component of $S$ on which $f$ is constant, and that
contains only two special points, contracts it.

The stratum $\delta$ has two preimages under the forgetful
map. The first one, 
$\widehat \delta_1$, consists of the functions defined on 
$2$-component curves, with the second pole and the distinguished
point on one component and the first pole on the other.
The second one, 
$\widehat \delta_2$ also consists of the functions defined on 
$2$-component curves, but with the second pole 
on one component and the first pole and the distinguished
point on the other. We also denote by $\widehat s(L)$
the preimage of $s(L)$.

The multiplicities of intersections of $\widehat \delta_1$
and $\widehat \delta_2$ with $\widehat s(L)$ are computed as
in the previous proposition. They are equal to, respectively,
$q(L) (a_1+1)$ and $q(L) (a_2+1)$.

The multiplicity of intersection of $\delta$ with $s(L)$ is
equal to their sum, that is, $(a_1+a_2+2) q(L)$.
\qed

\section{Proofs of the relations on generating functions}
\label{Sec:proofs}

Once we have studied the intersection between a primitive
or a simple stratum and $\Delta_{p,q}$ we can prove
Theorems~\ref{Thmmain2} and~\ref{Thmmain1}.

\paragraph{Proof of Theorem~\ref{Thmmain2}.}
Consider a monomial $Y$ in variables $t, y_1, y_2, \dots$.
Suppose its degree in $t$ equals $n$; thus it corresponds
to a stratum in $\P\cH_n$. The coefficient of this monomial
in $G$ is equal, up to a combinatorial coefficient, 
to the degree of the $LL$ map
of the corresponding stratum $\Sigma(Y)$.

Let us look at the ways in which the monomial $Y$ can arise
on the right-hand side of the equality of Theorem~\ref{Thmmain2}. 
Consider a decomposition list $L= (a_1, \dots, a_k; i,j)$
and $k$ monomials in the generating function $G$.
This data determines exactly one term on the
right-hand side of the equality of Theorem~\ref{Thmmain2}.
We suppose that this term contributes to the coefficient of
$X$. We denote by $p$ the degree in $t$ of the $k$th monomial.
It is easy to see that the same data also determines
a $k$-eared stratum that contains a component of intersection
of $\Sigma(Y)$ and $\Delta_{p,q}$.

Recall that the combinatorial coefficient, by which
the degree of $LL$ is multiplied in $G$, equals
$$
\frac1{n! \, |\Aut (Y)|}.
$$
As we multiply the monomials, the factors $n!$ are combined
into a multinomial coefficient
$$
n \choose {n_1, \dots, n_k}
$$
that enumerates the ways to distribute the $n$ numbered
poles among the $k$ components. Similarly, the coefficients
$|\Aut (Y)|$ enumerate the number of ways to distribute
the critical points of the same multiplicity among the
$k$ components. 

As we apply the operator $D_{a_i}$ to the $i$th
monomial, it is multiplied by its degree in the variable
$y_{a_i}$. This corresponds to the number of choices of
a critical point of multiplicity $a_i$ that will become
the intersection point of the $i$th peripheral component
with the central component in the $k$-eared stratum.

The coefficient $m(L)$ is the multiplicity of the intersection
of $\Sigma(Y)$ with $\Delta_{p,q}$ (see 
Propositions~\ref{Prop:deltacaps(L)k} and~\ref{Prop:deltacaps(L)2},
and Remark~\ref{Rem:2andk}).

Finally, according to the formula of 
Corollary~\ref{Cor:recurrence}, this should be multiplied by
$pq$ and divided by $2n(n-1)$. The multiplication by $pq$
is achieved by applying $D_0$ to the $k$th monomial and
to the product of the first $k-1$ monomials. Instead of
dividing the resulting monomial by $2n(n-1)$, we multiply
by $2n(n-1)$ the corresponding monomial on the left-hand
side of the equality, by writing
$$
2 \, t^2 \, \frac{\partial^2 G}{\partial t^2}
$$
instead of just $G$.

Thus each term on the right-hand side of the equality
corresponds to a $k$-eared stratum that contains
a component of $\Sigma(Y) \cap \Delta_{p,q}$, and correctly
computes the contribution of this component. 
All the components of the intersection are
taken into account. This finishes the proof.
\qed

\paragraph{Proof of Theorem~\ref{Thmmain1}.} 
Consider a monomial $X$ in variables $t, x_1, x_2, \dots$.
Suppose its degree in $t$ equals $n$; thus it corresponds
to a stratum in $\P\cH_n$. The coefficient of this monomial
in $F$ is equal, up to a combinatorial coefficient, 
to the degree of the $LL$ map
of the corresponding stratum $\Sigma(X)$.

Let us look at the ways, in which the monomial $X$ can arise
on the right-hand side of the equality of Theorem~\ref{Thmmain1}. 
The right-hand side of the equality of Theorem~\ref{Thmmain1}
contains the term
$$
\left( t \frac{\partial F}{\partial t} \right)^2
$$
and the terms
$$
(i+j) x_{i+j} 
\frac{\partial F}{\partial x_i}
\frac{\partial F}{\partial x_j}
$$
for all $i,j$. 

Let us choose one of the above terms and two monomials
of the function $F$. We denote by $p$ and $q$ their
degrees in the variable $t$. We suppose that as we apply the
chosen term to the two monomials we obtain a contribution
to the coefficient of the monomial $X$.

Then the choice of the term in the equation and of two
monomials determines a $2$-eared
stratum that constitutes a component of the intersection
$\Sigma(X) \cap \Delta_{p,q}$. The term 
$(t \, \partial F/ \partial t)^2$ corresponds to the case 
when two simple critical values get glued together. The terms of
the second type correspond to the case when the multiple
critical value is glued together with a simple critical value;
a cycle of length $i+j$ is then multiplied by a transposition
and is decomposed into a cycle of length $i$ and a cycle
of length $j$.

Recall that the combinatorial coefficient by which the
degree of $LL$ is multiplied to obtain the coefficient
of $F$ equals
$$
\frac1{(n-1)! \, |\Aut (X)|} = \frac{n}{n! \, |\Aut (X)|}.
$$
As we multiply two monomials, the factors $n!$ combine to
$$
{n! \choose p!\,  q!},
$$
which enumerates the ways to distribute the $n$ numbered
poles among the two components of the $2$-eared stratum.
Similarly, the factors $|\Aut (X)|$ account for the
number of ways to distribute the preimages of the only
multiple critical value among the two components. 

As we apply the operator $t \partial/ \partial t$ to
a monomial, we multiply it by the degree
of the stable functions in the corresponding stratum.
This gives the number of ways to choose a preimage
of a regular (noncritical) value of our function.
The chosen point will be the intersection point
between the two components of the nodal curve in
our $2$-eared stratum.
Similarly, when we differentiate a monomial with respect
to the variable $x_i$, this counts the number of ways
to choose a preimage of multiplicity $i$ of the
only multiple critical value. The chosen preimage will
again be the intersection point of the two components
of the nodal curve.

The coefficient $i+j$ is the multiplicity of the intersection
of $\Sigma(X)$ and $\Delta_{p,q}$ along the $2$-eared stratum
(Proposition~\ref{Prop:deltacaps(L)2}). According to the formula
in Corollary~\ref{Cor:recurrence}, it remains to multiply
the coefficient of our monomial by $pq$, and to divide
it by $2n(n-1)$. The factor $pq/n$ appears automatically
because of the factor $n$ in the numerator of the
combinatorial coefficient
$$
\frac{n}{n! \, |\Aut (X)|}
$$
by which the degree of $LL$ is multiplied. Instead of dividing
the result by $2(n-1)$, we multiply by $2(n-1)$ the
corresponding monomial on the left-hand side of
the equality, by writing 
$$
2t^2 \frac{\partial}{\partial t}
\left(\frac{F}{t} \right)
$$
instead of just $F$.

Thus the terms on the right-hand side enumerate all
the components of the intersection $\Sigma(X) \cap \Delta_{p,q}$
and each term computes correctly the contribution of
the corresponding component.

This finishes the proof. \qed

\section{Combinatorial results}
\label{y1y2}

\subsection{Some explicit formulas for Hurwitz numbers} 

Here we give some explicit formulas for the coefficients
of the function $G$. Some of them were known before,
but others are new.

We start by listing the already known expressions for
the coefficients of $G$ and the corresponding Hurwitz numbers.

In Theorem~\ref{Thm:6} below $y_*$ means ``any $y_i$'';
two asterisks do not necessarily denote the same subscript.
As usual, $|\Aut|$ is the number of ways to permute the variables
$y_i$ preserving their subscripts $i$.

\begin{theorem} \label{Thm:6}
\

1. The coefficient of 
$
\; \, t^n \underbrace{y_{n-1} y_* \dots y_*}_k \;
$
in $G$ equals
$
\; n^{k-3}/|\Aut|.
$
The corresponding Hurwitz number equals
$\; n^{k-3} \;$.

2. The coefficient of 
$
\; \, t^n \underbrace{y_{n-2} y_* \dots y_*}_k \;
$
in $G$ equals
$$
\frac{(k-2)\,(n-1)^{k-3}}{|\Aut|}.
$$
The corresponding Hurwitz number equals
$\; (k-2)\,(n-1)^{k-3}\;$.

3. The coefficient of 
$
\; \, t^n y_i y_1^{2n-2-i} \;
$
in $G$ equals, for $i \not= 1$,
$$
\frac{n^{n-i-3} \, (i+1)^{i+1}}{(n-i-1)! \, (i+1)!}.
$$
The corresponding Hurwitz number equals
$$
\frac{n^{n-i-3} \, (i+1)^{i+1} \, (2n-2-i)!}
{(n-i-1)! \, (i+1)!}.
$$

4. The coefficient of 
$
\; \, t^n y_1^{2n-2} \;
$
in $G$ equals
$
\; n^{n-3}/n! \; .
$
The corresponding Hurwitz number equals
$\; n^{n-3} \, (2n-2)!/n! \;$.
\end{theorem}

\paragraph{Proof.}\

1. The strata in questions are polynomial strata:
if we make a change of variables so that the critical
point corresponding to the factor $y_{n-1}$ is sent to 
$\infty$, all the functions of these strata will become
polynomials. The degrees of the Lyashko-Looijenga
map in the case of polynomials were found 
in~\cite{LanZvo}.

2. The strata in question are some particular strata
of the versal deformation of the singularity $D_n$. The
degree of the Lyashko-Looijenga map on these
strata was found in~\cite{PanZvo}.

3 and 4. These are primitive strata and the corresponding
degree is given by the Hurwitz formula (see
Section~\ref{primitive}).

\bigskip

Now we give some new formulas for the coefficients of
$G$ and the corresponding Hurwitz numbers.
We start with the following example. 

\begin{example}
Consider the stratum $\Sigma$ in $\P\cH_n$ such that a generic function
$f \in \Sigma$ has two double critical points while the remaining
$2n-6$ critical points are simple. The Hurwitz number corresponding
to the stratum $\Sigma$ equals
$$
\frac34 \, 
\frac{(2n-6)!}{n!} \cdot  n^{n-5} \cdot (n-1)\, (n-2)\,(27n^2-137n+180).
$$
\end{example}

This example has the following generalization.

Let $M$ be a monomial in variables $y_2$, $y_3$, \dots
(without $y_1$). Denote by $k$ the sum of the indices of
its factors and by $r$ the number of the factors (the
total degree of the monomial). Further denote by $s$
the maximum among the indices of the factors of $M$ and
the number $[(k+1)/2]$ (where $[\cdot]$ is the integer part).

\begin{conjecture} \label{Conj:pol}
The coefficient of 
$
\; \, t^n y_1^{2n-2-k} M \;
$
in $G$ equals
$$
\frac{n^{n-3-k+r}}{n!} \, P_M(n) \, ,
$$
where $P_M$ is a polynomial of degree $k$ with rational
coefficients divisible by $(n-1) \dots (n-s)$. The 
degree of the $LL$ map on the corresponding stratum equals
$$
|\Aut(M)| \;
\frac{(2n-2-k)!}{n!} \, n^{n-3-k+r}
\, P_M(n) \, .
$$
Here, as usual, $|\Aut(M)|$ is the number of ways to permute
the factors of $M$ preserving their indices.
\end{conjecture}

Although we do not know how to prove Conjecture~\ref{Conj:pol},
there is an algorithm that checks that the conjecture is
true for any given monomial $M$.

Here are some simplest polynomials $P_M$:
\begin{eqnarray*}
P_1(n) & = & 1, \\
P_{y_2}(n) & = & \frac92 \cdot (n-1)(n-2), \\
P_{y_3}(n) & = & \frac92 \cdot (n-1)(n-2)(n-3), \\
P_{y_4}(n) & = & \frac{32}3 \cdot (n-1)(n-2)(n-3)(n-4), \\
P_{y_2^2}(n) & = & \frac38 \cdot (n-1)(n-2)(27n^2-137n+180), \\
P_{y_2y_3}(n) & = & 8 \cdot (n-1)(n-2)(n-3)(6n^2-37n+60), \\
P_{y_3^2}(n) & = & \frac 29 \cdot
(n-1)(n-2)(n-3)(256n^3-2787n^2+10448n-13440), \\
P_{y_2^3}(n) & = & \frac1{48} \cdot 
(n-1)(n-2)(n-3)(729n^3 - 6723n^2+21026n-22680).
\end{eqnarray*}

We have not been able to find a general pattern for 
these polynomials.

\subsection{The algebra of generating functions}

Let us fix $c$ points on $\CP^1$ and assign a positive
integer $d_i$ and a partition $\kappa_i$ of $d_i$,
$1 \leq i \leq c$, to each marked point.
Consider $n$-sheeted coverings of $\CP^1$ by $\CP^1$
with simple preimages of $\infty$,
ramified over the $c$ distinguished points with
ramification types $\kappa_i$ and, in addition,
having $d(n) = 2n-2 - \sum d_i$ other fixed simple
ramification points.

Denote by $h_{\kappa_1, \dots, \kappa_c}(n)$ 
the corresponding Hurwitz number
(the number of the above coverings counted with
coefficient $1/|\Aut|$).

We regroup these Hurwitz numbers into a generating function
$$
f_{m_1, \dots, m_c}(t) = 
\sum_n \frac{h_{\kappa_1, \dots, \kappa_c}(n)}{d(n)!} t^n.
$$

We will prove that all such generating functions belong
to the following algebra $\A$.

\begin{definition}
Denote by $\A$ the subalgebra of the algebra of power
series in one variable $t$, generated by the series
$$
Y = \sum_{n \geq 1} \frac{n^{n-1}}{n!} \, t^n
\quad \mbox{and} \quad
Z = \sum_{n \geq 1} \frac{n^n}{n!} \, t^n
$$
\end{definition}

We start by a precise description of the algebra
$\A$.

\begin{proposition}
The algebra $\A$ is isomorphic to $\C[X, X^{-1}]$,
where $X=1-Y$, $X^{-1} = 1+Z$.
\end{proposition}

\paragraph{Proof.} All we have to prove is the equality
$(1-Y)(1+Z) = 1$. This equality follows from the identity
$$
\mathop{\sum_{p+q=n}}_{p,q \geq 1} \frac{n!}{p! \, q!} \, p^p \, q^{q-1}
=
(n-1)n^{n-1}, 
$$
which follows from Abel's identity
$$
\mathop{\sum_{p+q=n}}_{p,q \geq 0} 
\frac{n!}{p! \, q!} \, (a+p)^p \, (b+q)^{q-1}
= \frac{(a+b+n)^n}{b}, 
$$
by letting $a$ tend to $0$, then subtracting $n^n/b$ from both sides 
of the identity, and finally letting $b$ tend to $0$. 
For a proof of Abel's identity 
see~\cite{Abel}, \cite{Comtet}; see also \cite{EkhMaj}
for a computer-generated proof.
\qed

\bigskip

Denote by $A_n$ the following number
$$
A_n = \sum_{p+q=n} \frac{n!}{p! \, q!} \, p^p \, q^q \,
= n! \sum_{k = 0}^{n-2} \frac{n^k}{k!}. 
$$
For some combinatorial properties of this sequence
see~\cite{RioSlo}. (Exercise: using Abel's identity
prove that the two expressions defining $A_n$ are equal.)

\begin{proposition} \label{Prop:alg1}
The algebra $\A$ is formed by the formal power series of the form
$$
a + \sum_{n \geq 1} \frac{L(n) n^n + P(n) A_n}{n!} \, t^n,
$$
where $L$ is a Laurent polynomial and $P$ just a
polynomial in $n$.
\end{proposition}

\paragraph{Proof.} We will actually prove two more
precise statements.

(i) The subalgebra of $\A$ generated by $X$ (or by $Y$)
is the algebra of power series as above such that
the Laurent polynomial $L$ contains only negative
powers of $n$ and the polynomial $P$ is equal to $0$.

(ii) The subalgebra of $\A$ generated by $X^{-1}$ (or by $Z$)
is the algebra of power series as above such that
the Laurent polynomial $L$ is actually a polynomial.

Assertion (i) follows from the equality
$$
Y^k = k \sum_{n \geq 1} 
\frac{(n-k+1) \dots (n-1) \, n^{n-k}}{n!} \, t^k,
$$
which is proved by induction on $k$. 
Denote by $D$ the differential operator
$D = t\, \partial/ \partial t$.  It acts on power series by
multiplying the coefficient of $t^n$ by $n$.
The step of induction is carried out using the equality
$$
DY = Z = \frac{Y}{1-Y}.
$$
From this equality we deduce
$$
D \left( \frac{Y^{k+1}}{k+1} - \frac{Y^k}{k} \right) 
= (Y^k-Y^{k-1}) \, \frac{Y}{1-Y} =-Y^k.
$$
Knowing that the free term of $Y^{k+1}$ equals $0$,
we can deduce the expression of $Y^{k+1}$ from that
of $Y^k$.

The proof of assertion (ii) is similar, although, 
as far as we know, there is no simple formula for $Z^k$. 

By definition of the sequence $A_n$, we have
$$
\sum_{n \geq 1} \frac{A_n}{n!} \, t^n = Z^2.
$$
Therefore what we must prove is that the series
$D^kZ$ and $D^k(Z^2)$ for $k \geq 0$
can be expressed as polynomials in $Z$ and vice versa.

From $(1-Y)(1+Z) = 1$ and $Z = DY$ we find
$$
DZ = Z(1+Z)^2.
$$
It follows by induction that $D^kZ$ can be expressed
as a polynomial in $Z$ of degree $2k+1$ with positive
integer coefficients. Hence $D^k (Z^2)$ can be expressed as a 
polynomial in $Z$ of degree $2k+2$ with positive integer 
coefficients. Since we have obtained exactly once every
possible degree of polynomials in $Z$, it follows that
any polynomial in $Z$ is a linear combination of the
series $D^kZ$ and $D^k (Z^2)$.
\qed

We do not know whether the algebra $\A$ has already
appeared in the literature. It certainly plays an
important role in the understanding of ramified coverings
of the sphere, as shows Theorem~\ref{Thm:algebra} from
Introduction. We restate it here.

\setcounter{theorem}{1}

\begin{theorem}
For any partitions $\kappa_1, \dots, \kappa_c$, the series 
$f_{\kappa_1, \dots, \kappa_c}(t)$
belongs to the algebra $\A$. 
\end{theorem}

\begin{conjecture}
If we write the series $f_{\kappa_1, \dots, \kappa_c}(t)$
in the form of Proposition~\ref{Prop:alg1}, then the corresponding
polynomial $P(n)$ is equal to $0$.
\end{conjecture}

Why the polynomial $P$ always equals $0$ remains mysterious.
In the proof of Theorem~\ref{Thm:algebra} the function
$f_{\kappa_1, \dots, \kappa_c}$ is obtained as a sum
of products of similar functions. Another way to obtain
$f_{\kappa_1, \dots, \kappa_c}$ as a different sum of
products of similar functions can be obtained from
a generalization of Theorem~\ref{Thmmain2}.
(The relation of Theorem~\ref{Thmmain2} itself gives
a method for calculating the elements of the
algebra $\A$ that correspond to various generating functions
$f_{\kappa_1, \dots, \kappa_c}(t)$ corresponding to
simple strata.) However, in both methods, the polynomials
$P$ that correspond to various terms of the sum
are not equal to $0$. They cancel out only after the
addition.

\begin{remark}
Theorem~\ref{Thm:algebra} resembles very much another theorem
on the ramified coverings of a torus.
Fix a finite number of ramification points
on a torus and consider $n$-sheeted coverings of
the torus having ramifications of prescribed types
over the marked points. The numbers of such ramifications
can be regrouped in a generating function. It is known
(see~\cite{BloOko}) that this generating function will always
belong to a very particular subalgebra of the algebra
of power series, namely, to the algebra of quasimodular
forms (that is, the algebra generated by the three
Eisenstein series $E_2$, $E_4$, $E_6$). This important
result was predicted by mirror symmetry for elliptic curves.
\end{remark}

\paragraph{Proof of Theorem~\ref{Thm:algebra}.}
We will prove the theorem by induction on $c$. First
of all note that the base of induction (the cases $c=0,1$)
is an immediate consequence of the Hurwitz formula
(see~\cite{Hurwitz2}):
$$
h_\kappa=\frac{(2n-2-d(\kappa))!}{(n-d(\kappa)-m)! \, |\Aut(\kappa)|}
\prod_{i=1}^m\frac{(k_i+1)^{k_i+1}}{(k_i+1)!}n^{n-d(\kappa)-3},
$$
where $\kappa = (k_1,\dots,k_m)$, $d(\kappa) = k_1 + \dots + k_m$,
and $|\Aut(\kappa)|$ is the number of ways to permute the
numbers $k_i$ preserving their values.

Now we will show that each generating
function $f_{\kappa_1, \dots, \kappa_c}(t)$ can
be obtained from similar generating functions involving
$c-1$ instead of $c$ partitions by applying sums, products, 
and the operator
$D = t \, \partial/\partial t$. Since all these operations
preserve the algebra $\A$, if all the generating functions
for $c-1$ partitions belong to $\A$, so does the final 
generating function. The proof goes in the spirit of~\cite{GouJac}.

The monodromy of a ramification point with ramification
type $\kappa = (k_1, \dots, k_m)$ is a permutation whose
cycle lengths are $k_1+1, \dots, k_m+1$ and, in addition,
some number of cycles of length~$1$. For simplicity we
will call $\kappa$ the {\em cycle type} of the permutation.
Consider a factorization of the identity permutation on 
$n$ elements into a product of $c$ permutations 
$\sigma_1, \dots, \sigma_c$ with cycle types
$\kappa_1, \dots, \kappa_c$ and $d(n)$ transpositions
$\tau_1, \dots, \tau_{d(n)}$ (in that order)
such that the group generated by all these permutations
acts transitively on the $n$ elements. Such a factorization
is called a {\em transitive factorization}. It is easy
to see that the number $d(n)$ is the smallest possible
number of transpositions that should be added to
the other permutations so that
a transitive factorization would exist. It follows that every
transitive factorization is the set of monodromies of
a ramified covering of a sphere by a sphere.
Moreover, the number of all transitive factorizations equals
$n! \, h_{\kappa_1, \dots, \kappa_c}(n)$. 

Before proceeding with the proof, let us give an example
of what will happen. Consider the generating function
 $f_{1,1}$ corresponding to the case when $c=2$ 
and the two special critical values
are actually just simple critical values; in other
words, the corresponding monodromies are transpositions.
A product of two transpositions can be either a cycle
of length~3 (corresponding to a double critical point),
or a permutation with two cycles of lengths~2 (a critical
value attained at two simple critical points), or
the identity permutation with two distinguished elements
(the ramified covering splits into two irreducible
components that
intersect at one point). Therefore we will need three
simpler generating functions: $f_2$, $f_{1^2}$, and $f_{\emptyset}$
(the latter corresponds to the case $c=0$). We have 
the following equality:
$$
f_{1,1} = 3f_2 + 2 f_{1^2}+ \frac12(D f_{\emptyset})^2.
$$
The coefficient $3$ is due to the fact that there are {\em three}
ways do decompose a cycle of length $3$ into a product
of $2$ transpositions. Similarly, the coefficient $2$ is due to
the fact that a permutation with $2$ cycles of length~$2$
can be decomposed into a product of $2$ transpositions
in {\em two} ways. The operator $D$ is applied to the function
$f_{\emptyset}$ because if there are
$p$ and $q$ sheets in the two transitivity components,
then there are, respectively, $p$ and $q$ ways to choose
their intersection point on each of the components
(and the operator $D$ multiplies by $p$ the coefficient
of $t^p$). Finally, the coefficient $1/2$ is due to the
fact that the two transitivity components are indistinguishable
from each other. The reader can check that the equality
is true using the expressions

\begin{eqnarray*}
f_{\emptyset} & = & \sum_{n \geq 1} \frac{n^{n-3}}{n!} t^n 
=  -\frac1{12}(X-1)(2X^2+5X+5), \\
D f_{\emptyset} & = & \sum_{n \geq 1} \frac{n^{n-2}}{n!} t^n
=  -\frac12(X-1)(X+1), \\
f_{1^2} & = & 2 \sum_{n \geq 1} (n-3)(n-2)(n-1) \frac{n^{n-4}}{n!} t^n
=  \frac12 (X-1)^4, \\
f_2 & = & \frac92 \sum_{n \geq 1} (n-2)(n-1) \frac{n^{n-4}}{n!} t^n
=  -\frac18 (X-1)^3(3X+1), \\
f_{1,1} & = & \sum_{n \geq 1} (2n-2)(2n-3) \frac{n^{n-3}}{n!} t^n
=  -\frac12 (X-1)^2(2X-3).
\end{eqnarray*}

We will show that similar relations hold in the general
case. We do not need to specify the coefficients in these
relations; it suffices to know that they are some fixed
rational numbers. Now we proceed with the proof.

Let us see what happens to a transitive factorization
as we multiply the permutations $\sigma_{c-1}$ and $\sigma_c$.
We obtain a new list of permutations, whose
product is equal to the identity permutation. However
the group generated by these permutations is not
necessarily transitive. Denote by $k$ the number of
its transitivity components. Each transitivity component
is a subset of the set $\{ 1, \dots, n \}$.
For each component we actually obtain a transitive factorization
of the identity permutation acting on this subset. The elements
of each of the partitions $\kappa_1, \dots, \kappa_{c-2}$
are distributed among the transitivity components. The $d(n)$
transpositions are also distributed among the transitivity
components. As for the product $\sigma_{c-1} \sigma_c$, we
will separate its cycles into two groups. The first group
consists of those cycles that have at least one element in common
with a cycle of length greater than~$1$ either
in $\sigma_{c-1}$ or in $\sigma_c$. The cycles in this
first group can have length~$1$ (singletons) or greater than~$1$.
The cycles of length~$1$ will be called {\em special} singletons.
The second group consists of all the other cycles. All these
cycles are automatically singletons. If the number $n$
is big enough, the second group will be much bigger than
the first one. The cycles of length greater than~$1$ and
the special singletons of the product
$\sigma_{c-1} \sigma_c$ are again distributed in some
way among the $k$ transitivity components. Thus to every
transitive factorization with cycle types 
$\kappa_1, \dots, \kappa_c$ we can assign unambiguously the
following data: $k$ transitive factorizations acting on some subsets of
the set $\{1, \dots, n\}$; each of these transitive
factorizations is composed of $c-1$ permutations with
some given cycle types, plus the necessary number of
transpositions; some number of special singletons is
specified among the singletons of the $(c-1)$th permutation.

Working backwards, one can see that the above data allows
us to reconstruct several possible transitive factorizations
with cycle types $\kappa_1, \dots, \kappa_c$. The number of
possibilities is equal to the product of two factors. 

The first factor is a multinomial
coefficient that counts the number of ways to order the
$d(n)$ transpositions. These transpositions are already
ordered inside each of the $k$ transitive factorizations,
and we have to order them all respecting these partial
orders. As we multiply two generating functions, this
first factor appears automatically, because of the
coefficient $d(n)!$ in the denominator of the
coefficient of $t^n$ in $f_{\kappa_1, \dots, \kappa_c}$.

The second factor is the number of ways to reconstruct the
permutations $\sigma_{c-1}$ and $\sigma_c$ knowing their
product. More precisely, the $(c-1)$th permutations
in the $k$ transitive factorizations can be regrouped
into a unique permutation acting on the set $\{1, \dots, n\}$.
This permutation has some cycles of lengths greater than~$1$
and, in addition, a given number of special singletons. We
must decompose this permutation into a product of two
permutations $\sigma_{c-1}$ and $\sigma_c$ (with given
cycle types) in such a way that all nonspecial singletons
remain singletons in $\sigma_{c-1}$ and $\sigma_c$. Moreover,
the $k$ subsets of the subdivision must merge into a unique
transitivity component after we have replaced the product
$\sigma_{c-1} \sigma_c$ by two permutations $\sigma_{c-1}$
and $\sigma_c$. The number of such factorizations is the second
factor. This factor, divided by some product of factorials
due to indistinguishable transitivity components, will play
the role of the coefficient in the expression for
$f_{\kappa_1, \dots, \kappa_c}$.

Now we can precisely describe how the generating function
$f_{\kappa_1, \dots, \kappa_c}$ is expressed from similar
generating functions assigned to lists of $c-1$ partitions.
First we make a list of all possible ways to divide the
elements of each of the partitions $\kappa_1, \dots \kappa_{c-2}$ 
into $k$ parts and to distribute among the same $k$ parts
the cycle lengths and the special singletons of each
possible product of the permutations $\sigma_{c-1}$ and $\sigma_c$.
This gives us a finite list of possibilities; each possibility
can be described by a table:
$$
\begin{array}{cccc}
\kappa_1^{(1)}, & \dots, & \kappa_{c-1}^{(1)}, & m^{(1)}; \\
\kappa_1^{(2)}, & \dots, & \kappa_{c-1}^{(2)}, & m^{(2)}; \\
&\dots&&\\
\kappa_1^{(k)}, & \dots, & \kappa_{c-1}^{(k)}, & m^{(k)}.
\end{array}
$$
Here $\kappa_i^{(j)}$ for $i \leq c-2$
is the partition formed by those elements
of the partition $\kappa_i$ that were assigned to the
$j$th part (empty partitions are allowed); 
$\kappa_{c-1}^{(j)}$ is the similar partition assigned
to the restriction to the $j$th part of the
product $\sigma_{c-1} \sigma_c$; finally, $m^{(j)}$ is the
number of special singletons in the $j$th part.

To each row of this table we assign the generating function
corresponding to the partitions of this row. The
coefficients of this generating function count the
number of transitive factorizations in the corresponding
transitivity component. Further,
denote by $d^{(j)}$ the sum of the elements
of $\kappa_{c-1}^{(j)}$. Then we apply to the corresponding
generating function the operator 
$$
(D-d^{(j)}) \dots (D-d^{(j)}-m^{(j)}+1).
$$
This corresponds to choosing $m^{(j)}$ special singletons
among all the available singletons. Now we multiply the
obtained expressions for each row. The product
thus obtained is further multiplied by the
number of ways to reconstruct the permutations
$\sigma_{c-1}$ and $\sigma_c$ knowing the
permutation $\sigma_{c-1}\sigma_c$, and then 
divided by the number of automorphisms of the set of rows
of the table. The sum of such terms over all possible
tables gives us an expression for the generating
function $f_{\kappa_1, \dots,\kappa_c}$.

This finishes the proof. \qed

\end{document}